\documentclass{amsart}
\title{Borel's Conjecture in topological groups} 
\author{Fred Galvin and Marion Scheepers}
\date{}

\usepackage{amsfonts}
\usepackage{amsmath, amsthm}
\usepackage{colonequals}
\usepackage{txfonts}
\newtheorem{theorem}{{\bf Theorem}}
\newtheorem{proposition}[theorem]{{\bf Proposition}}
\newtheorem{lemma}[theorem]{{\bf Lemma}}
\newtheorem{corollary}[theorem]{{\bf Corollary}}

\newtheorem{problem}{{\bf Problem}}

\newcommand{\naturals}{{\mathbb N}}

\newcommand{\reals}{{\mathbb R}}

\newcommand{\identity}{{\mathbf 1}}
\newcommand{\forces}{\mathrel{\|}\joinrel\mathrel{-}}
\newcommand{\changarrow}{\twoheadrightarrow}
\newcommand{\poset}{{\mathbb P}}
\newcommand{\sone}{{\sf S}_1}

\newcommand{\open}{\mathcal{O}}
\newcommand{\onbd}{\open_{nbd}}

\makeatletter

\begin{document}
\maketitle
\begin{abstract}
We introduce a natural generalization of Borel's Conjecture. For each infinite cardinal number $\kappa$, let {\sf BC}$_{\kappa}$ denote this generalization. Then ${\sf BC}_{\aleph_0}$ is equivalent to the classical Borel conjecture. Assuming the classical Borel conjecture, $\neg{\sf BC}_{\aleph_1}$ is equivalent to the existence of a Kurepa tree of height $\aleph_1$. Using the
connection of ${\sf BC}_{\kappa}$ with a generalization of Kurepa's Hypothesis, we obtain the following consistency results:
\begin{enumerate}
  \item{If it is consistent that there is a 1-inaccessible cardinal then it is consistent that ${\sf BC}_{\aleph_1}$.}
  \item{If it is consistent that ${\sf BC}_{\aleph_1}$, then it is consistent that there is an inaccessible cardinal.}
  \item{If it is consistent that there is a 1-inaccessible cardinal with $\omega$ inaccessible cardinals above it, then $\neg{\sf BC}_{\aleph_{\omega}} \, +\, (\forall n<\omega){\sf BC}_{\aleph_n}$ is consistent.}
  \item{If it is consistent that there is a 2-huge cardinal, then it is consistent that ${\sf BC}_{\aleph_{\omega}}$.} 
  \item{If it is consistent that there is a 3-huge cardinal, then it is consistent that ${\sf BC}_{\kappa}$  for a proper class of cardinals $\kappa$ of countable cofinality.} 
\end{enumerate}  
 \end{abstract}

A metric space $(X,d)$ is \emph{strong measure zero} if there is for each sequence $(\epsilon_n:n<\omega)$ of positive real numbers a corresponding sequence $(U_n:n<\omega)$ of open sets such that for each $n$ the set $U_n$ has $d$-diameter at most $\epsilon_n$, and $\{U_n:n<\omega\}$ covers $X$. Strong measure zero metric spaces are necessarily separable. E. Borel \cite{Borel} conjectured that strong measure zero sets of real numbers are countable. The metric notion of strong measure zero has a natural generalization to non-metric contexts. Rothberger \cite{rothberger38} introduced a generalization to the class of topological spaces. We consider a generalization to the class of topological groups. Most of our results can be presented in the more general context of uniformizable spaces, but we found no advantage to presenting it thus. 

Borel's Conjecture also has natural generalizations to these non-metric contexts. These generalizations expose, as in the metric case, interesting connections with the foundations of mathematics. The generalization of Borel's Conjecture considered here is quite different from what Halko and Shelah considered in \cite{HS}.

After a brief introduction of notation and terminology we define \emph{Rothberger boundedness} and \emph{Rothberger spaces} in Section 1. In Section 2 we introduce a generalization of Borel's Conjecture and in Section 3 explore connections between it and other combinatorial structures. In Section 4 we give a number of consistency and independence results regarding the generalization introduced in the earlier sections. 

By a well-known theorem of Kakutani a topological group is ${\sf T}_0$ if, and only if, it is ${\sf T}_{3\frac{1}{2}}$\footnote{A ${\sf T}_0$ topological group need not be normal.}. Throughout this paper we shall assume, without further notice, that all groups considered are ${\sf T}_{3\frac{1}{2}}$. Correspondingly, all topological spaces we consider here are assumed to be ${\sf T}_{3\frac{1}{2}}$.

Let $(G,*)$ be a topological group with identity element $\identity$. 
For nonempty subsets $A$ and $B$ of $G$ and for $g\in G$ the symbol $A*B$ denotes the set $\{a*b:a\in A\mbox{ and }b\in B\}$, and $g*A$ denotes $\{g\}*A$. The symbol $\mathcal{O}$ denotes the set of all nonempty open covers of $G$. 

Let $U$ be an open neighborhood of $\identity$. Then $\open(U) = \{g*U:g\in G\}$ is an open cover of $G$. Define: 
\[
  \open_{nbd}:=\{\open(U): U \mbox{ an open neighborhood of } \identity\}.
\]
According to Guran \cite{Guran} the topological group $(G,*)$ is $\aleph_0$-\emph{bounded} if each element of $\onbd$ has a countable subset which covers $G$. A topological group is said to be \emph{pre-compact} if each element of $\open_{nbd}$ has a finite subset covering the group.

While pre-compact topological groups and Lindel\"of topological groups are $\aleph_0$-bounded the converse is not true. 
The class of $\aleph_0$-bounded groups has nice preservation properties: Every subgroup of an $\aleph_0$-bounded group is $\aleph_0$-bounded, any (finite or infinite) product of $\aleph_0$-bounded groups is $\aleph_0$-bounded, every continuous homomorphic image of an $\aleph_0$-bounded group is $\aleph_0$-bounded, and if a dense subgroup of a group is $\aleph_0$-bounded, then so is the group. The survey \cite{MT} gives a good introduction to $\aleph_0$-bounded groups.
\begin{theorem}[Guran]\label{guranembth} A topological group is $\aleph_0$-bounded if, and only if, it embeds as a topological group into a product of second countable topological groups.
\end{theorem}
By the Birkhoff-Kakutani theorem each second countable ${\sf T}_0$ topological group is metrizable. Thus the groups appearing as factors in the product in Guran's theorem are separable and metrizable. 
Guran's theorem has the following quantified form (see \cite{MT}):
\begin{theorem}\label{guranquantify} For an $\aleph_0$-bounded topological group $(G,*)$ and an infinite cardinal number $\kappa$ the following are equivalent:
\begin{enumerate}
  \item{The weight of $G$ is $\kappa$.}
  \item{The character of $G$ is $\kappa$.} 
  \item{$\kappa$ is the smallest infinite cardinal such that $G$ embeds as a topological group into a product of $\kappa$ separable metrizable topological groups.}
\end{enumerate}
\end{theorem}

\section{Rothberger boundedness in topological groups, Rothberger spaces.}

For collections $\mathcal{A}$ and $\mathcal{B}$ 
the symbol $\sone(\mathcal{A},\mathcal{B})$ denotes the selection principle 
\begin{quote} For each sequence $(A_n:n\in\naturals)$ of elements of $\mathcal{A}$ there is a sequence $(b_n:n\in\naturals)$ such that for each $n$, $b_n\in A_n$, and $\{b_n:n\in\naturals\}\in\mathcal{B}$.
\end{quote}
A topological space is said to be a \emph{Rothberger space} if it has the property $\sone(\open,\open)$ (for an introdction to Rothberger spaces the reader could consult \cite{coc2}). The topological group $(G,*)$ is said to be \emph{Rothberger bounded} if it has the property $\sone(\open_{nbd},\open)$. For a \emph{subset} $X$ of the topological group $(G,*)$, $\open_X$ denotes the family of covers of $X$ by sets open in $G$. $X$ is said to be Rothberger bounded if $\sone(\onbd,\open_X)$ holds. These concepts are named after Rothberger who introduced $\sone(\open,\open)$ and who considered a close analogue of this boundedness property in the \emph{Hilfssatz} on page 51 of his paper \cite{rothberger38}.  

If a subset of a topological group is a Rothberger space, then it is Rothberger bounded in the group. The converse is not true. A subspace of a Rothberger space need not be a Rothberger space, but subsets of Rothberger bounded sets are Rothberger bounded.
Rothberger boundedness of a subset of a group is preserved by continuous group homomorphisms and countable unions. The property of being a Rothberger space is preserved by continuous surjections and by countable unions. 
\begin{proposition}\label{0dim} Every Rothberger bounded subset of an $\aleph_0$-bounded topological group is zero-dimensional.
\end{proposition}
{\flushleft{\bf Proof:}} Let $(G,*)$ be an $\aleph_0$-bounded group. Choose by Guran's theorem separable metrizable groups $(G_i,*_i)$, $i\in I$ such that $(G,*)$ embeds as topological group in the product $\Pi_{i\in I}(G_i,*_i)$ and let $\Phi$ be an embedding. The projection of $\Phi\lbrack G\rbrack$ on each coordinate of this product is a metrizable group. The Rothberger boundedness of subsets of $G$ is also preserved by the composition of $\Phi$ and projections.

In metrizable groups Rothberger boundedness coincides with strong measure zero.  Thus a Rothberger bounded subset $X$ of an $\aleph_0$-bounded ${\sf T}_0$ group projects onto a metrizable strong measure set in each coordinate projection. By a theorem of Marczewski \cite{SM} strong measure zero metrizable spaces are zero-dimensional. Thus, as products and subspaces of zero-dimensional spaces are zero-dimensional, $X$ is zero-dimensional. $\Box$

The proof of Proposition \ref{0dim} shows: Borel's Conjecture implies that each Rothberger bounded subset of an $\aleph_0$-bounded  topological group embeds as a Rothberger bounded subset into a product of countable topological groups (see \cite{CarlsonBC} Theorem 3.2). It also follows that Rothberger spaces must be zero-dimensional, and that Borel's Conjecture implies that  Rothberger spaces embed into products of countable topological groups.

If $C\subseteq \kappa$ is nonempty and if $S$ is a subset of $\prod_{\alpha<\kappa}G_{\alpha}$, then $S_C =\{f\lceil_C:f\in S\}$.

\begin{lemma}\label{grouprestrictlemma} Let $\kappa$ be an infinite cardinal number. Let $(G_{\alpha}:\alpha<\kappa)$ be topological groups and let a subset $X$ of $G = \prod_{\alpha<\kappa}G_{\alpha}$ be given. The following are equivalent:
\begin{enumerate}
  \item{$X$ is Rothberger bounded.}
  \item{For each countable set $C\subseteq \kappa$ the set $X_{C}$ is a Rothberger bounded subset
        of $G_{C}$.}
\end{enumerate}
\end{lemma}
{\flushleft{\bf Proof:}} \underline{(1)$\Rightarrow$(2):}  A continuous group homomorphism preserves Rothberger boundedness.

\underline{(2)$\Rightarrow$(1):} Suppose for each countable $C\subseteq \kappa$ we have that $X_C$ is a Rothberger bounded subset of $G_{C}$. Let $(U_n:n<\omega)$ be a sequence of neighborhoods of the identity element of $G$. We may assume that each $U_n$ is a basic open set such that there is a finite set $F_n\subseteq \kappa$ and for each $x\in F_n$ a neighborhood $N_x$ of $\identity_x$ such that $U_n=\{f\in\,G:(\forall x\in F_n)(f(x)\in N_x)\}$. Let $C$ be a countably infinite subset of $\kappa$ for which $\bigcup_{n<\omega}F_n\subseteq C$. Then for each $n$, $V_n=U_n\lceil _C\subseteq\,G_{C}$ is a neighborhood of the identity element of $G_{C}$. Since $X_C$ is Rothberger bounded choose for each $n$ an $x_n\in\, G_{C}$ such that $X_C \subseteq \bigcup_{n<\omega} x_n* V_n$. For each $n$ choose $y_n\in\,G$ such that $y_n\lceil_C=x_n$. Then we have $X\subseteq \bigcup_{n<\omega}y_n* U_n$. It follows that $X$ is Rothberger bounded.
 $\Box$

In the case of Rothberger spaces Lemma \ref{grouprestrictlemma} has the following analogue:
\begin{lemma}\label{spacerestrictlemma} Let $\kappa$ be an infinite cardinal number. Let $(G_{\alpha}:\alpha<\kappa)$ be topological groups and let a subspace $X$ of $G = \prod_{\alpha<\kappa}G_{\alpha}$ be given. Then the following are equivalent:
\begin{enumerate}
  \item{$X$ is a Rothberger space.}
  \item{$X$ is Lindel\"of and for each countable set $C\subseteq \kappa$ the set $X_{C}$ is a Rothberger subspace 
        of $G_{C}$.}
\end{enumerate}
\end{lemma}

A Rothberger bounded $\sigma$-compact subset of a topological group is a Rothberger space: First note that a union of countably many Rothberger subspaces of a space is again a Rothberger subspace. Thus, it suffices to show that a compact Rothberger bounded subset of a topological group is a Rothberger space. For this, use of the following analogue of the Lebesgue covering Lemma, following from Theorem 6.33 in \cite{Kelley}:
 \begin{lemma}\label{lebesguecov} Let $C$ be a compact subset of a topological group $(G,*)$ and let $\mathcal{U}$ be a cover of $C$ by sets open in $G$. Then there is a neighborhood $N$ of the identity of $G$ such that for each $x\in C$ there is a $U\in\mathcal{U}$ such that $x*N\subseteq U$.
  \end{lemma}

\begin{corollary}\label{galvin1} For each infinite cardinal $\kappa$ any topological group $\prod_{\alpha<\kappa}G_{\alpha}$, where each $G_{\alpha}$ has at least two elements, has a Rothberger subgroup\footnote{That is, a subgroup which is a Rothberger space.}
 of cardinality $\kappa$.
\end{corollary}
{\flushleft{\bf Proof:}} 
For each $\alpha < \kappa$, choose a nontrivial (meaning that it has more than one element) countable subgroup $H_{\alpha}$ of $G_{\alpha}$, and let $H = \prod_{\alpha < \kappa}H_{\alpha}$.
Then $X = \{x \in H: x_{\alpha} = {\mathbf 1} \mbox{ for all but finitely many } \alpha\}$ is a subgroup of H of cardinality $\kappa$. Inasmuch as $X$ is Lindel\"of (in fact $\sigma$-compact), and $X_C$ is countable whenever $C$ is countable, it follows by Lemma 5 that X is a Rothberger space.
$\Box$

Thus there are Rothberger (and so Rothberger bounded) groups of all infinite cardinalities.

\section{The Generalized Borel Conjecture}

Let {\sf BC} denote Borel's conjecture that each strong measure zero set of real numbers is countable\footnote{Sierpi\'nski \cite{sierpinski} proved that the Continuum Hypothesis implies the negation of Borel's conjecture, and Laver \cite{laver} proved that Borel's conjecture is consistent relative to the consistency of {\sf ZFC}, the Zermelo-Fraenkel axioms plus the Axiom of Choice. Thus, Borel's conjecture is independent of {\sf ZFC}.}. For the real line with the addition operation, Borel's conjecture translates to the statement that every Rothberger bounded subset of the topological group $(\reals,+)$ is countable.

\begin{theorem}\label{BCequivs} The following statements are equivalent:
\begin{enumerate}
  \item{{\sf BC}}
  \item{Each strong measure zero metric space is countable.}
  \item{Each Rothberger bounded subset of a group of countable weight is countable.}    
  \item{Each subgroup\footnote{(4) remains equivalent to {\sf BC} if we change "subgroup" to "subfield".} of the real line, all of whose finite powers are Rothberger, is countable.}
\end{enumerate}  
\end{theorem}
{\flushleft{\bf Proof:}} (1)$\Leftrightarrow$(2): This result is due to T.J. Carlson \cite{CarlsonBC}.\\
(2)$\Rightarrow$(3): By the Kakutani-Birkhoff Theorem a ${\sf T}_0$ group of countable weight is metrizable by a left-invariant metric. Any Rothberger bounded subset $X$ of such a group is strong measure zero in such a left-invariant metric.
By 2), $X$ is countable.\\
(3)$\Rightarrow$(4): Consider a subgroup $G$ of the real line such that each finite power of $G$ is Rothberger. Since $G$ is Rothberger it is Rothberger bounded. Since the real line is a second countable group, 3) implies that $G$ is countable.\\
(4)$\Rightarrow$(1): 
If $X$ is a set of real numbers, then $\lbrack X\rbrack$, the subgroup of the real line generated by $X$, is a countable union of continuous images of finite powers of $X$, and the same goes for any finite power of $\lbrack X\rbrack$. Hence, if $X$ has the property that all of its finite powers are Rothberger, then $\lbrack X\rbrack$ also has that property. Thus the implication (4)$\Rightarrow$ (1) follows from the result of Tsaban and Weiss \cite{TW} that {\sf BC} is equivalent to the statement that each subset of the real line, all of whose finite powers are Rothberger, is countable.
$\Box$

For $\lambda$ a cardinal number and $(G,*)$ a topological group, ${\sf BC}(G,<\lambda$) states:
\begin{quote}
{\tt  Each Rothberger bounded subset of $(G,*)$ has cardinality less than $\lambda$}.
\end{quote}
${\sf BC}(G,<\lambda^+)$ is also written ${\sf BC}(G,\lambda)$ and ${\sf BC}(G,\omega)$ is also written ${\sf BC}(G)$. 

It is clear that if $\mu<\lambda$ then ${\sf BC}(G,<\mu)$ implies ${\sf BC}(G,<\lambda)$.
Moreover, if $H$ is a subgroup of the group $G$ then ${\sf BC}(G,<\mu)$ implies ${\sf BC}(H,<\mu)$. 

Corollary \ref{galvin1} shows that for each infinite cardinal $\kappa$ the statement ${\sf BC}(^{\kappa}2,<\kappa)$ is false. The status of ${\sf BC}(G,\kappa)$ for $\aleph_0$-bounded groups of weight $\kappa$ is more elusive. 
For an infinite cardinal number $\kappa$ we define, inspired by (3) of Theorem \ref{BCequivs}:
\vspace{0.1in}

\begin{tabular}{lcl}
${\sf BC}_{\kappa}$: & & {\tt Each Rothberger bounded subset of an $\aleph_0$-bounded group of} \\
                     & & {\tt weight $\kappa$ has cardinality at most $\kappa$}.
\end{tabular} 
\vspace{0.1in}

Thus, ${\sf BC}_{\aleph_0}$ is Borel's Conjecture, ${\sf BC}$. It is also clear that for each infinite cardinal $\kappa$, ${\sf BC}_{\kappa}$ implies ${\sf BC}(^{\kappa}2, \kappa)$. The status of ${\sf BC}(^{\kappa}2, \kappa)$ is the focus of this paper.

\section{${\sf BC}_{\kappa}$ for uncountable $\kappa$.}

Let $\lambda\le\kappa$ be uncountable cardinal numbers. A family $\mathcal{F}$ of subsets of $\kappa$ such that $\vert \mathcal{F}\vert >\kappa$ and for each infinite subset $A$ of $\kappa$ such that $\vert A\vert<\lambda$, we have $\vert\{X\cap A:X\in\mathcal{F}\}\vert\le \vert A\vert$, is said to be a $(\kappa,\lambda)$ \emph{Kurepa family}\footnote{This definition is like the one in Chapter VII.3 of \cite{KDevlin}, but we do not require $\kappa$ to be regular.}. The $(\kappa,\lambda)$ \emph{Kurepa Hypothesis}, {\sf KH}($\kappa$,$\lambda$), is the assertion that there exists a $(\kappa,\lambda)$ Kurepa family. ${\sf KH}(\aleph_1,\aleph_1)$ is the classical Kurepa Hypothesis.

\begin{theorem}\label{first} Let $\kappa$ be an uncountable cardinal.
Let $(G_{\alpha}:\alpha<\kappa)$ be a family of topological groups, each with more than one element. If ${\sf KH}(\kappa,\aleph_1)$,  then $\prod_{\alpha<\kappa}G_{\alpha}$ has a Rothberger bounded subset (indeed, subgroup) of cardinality $\kappa^+$. 
\end{theorem}
{\flushleft{\bf Proof:}} For each $\alpha<\kappa$ let $\identity_{\alpha}$ denote the identity element of, and let $g_{\alpha}$ be any other element of the group $G_{\alpha}$. Let $\mathcal{F}$ be a $(\kappa,\aleph_1)$ Kurepa family on $\kappa$. For each $X\in\mathcal{F}$ define $\phi_X\in\prod_{\alpha<\kappa}G_{\alpha}$ so that for each $\alpha<\kappa$
\[
  \phi_X(\alpha) = \left\{
                          \begin{tabular}{ll}
                           $\identity_{\alpha}$ & \mbox{ if $\alpha\not\in X$}\\
                           $g_{\alpha}$         & \mbox{ if $\alpha\in X$}\\
                          \end{tabular}
                   \right.
\]
Then $S =\{\phi_X:X\in\mathcal{F}\}$ is a subset of cardinality $\vert\mathcal{F}\vert$ of $\prod_{\alpha<\kappa}G_{\alpha}$.

For each countable subset $C$ of $\kappa$, the set $S_C =\{\phi_X\lceil_C:X\in \mathcal{F}\}$ has the same cardinality as $\{X\cap C:\, X\in\mathcal{F}\}$ and thus is countable. But then for each countable set $C\subset \kappa$, $S_C$ is a Rothberger, and thus Rothberger bounded, subset of $\prod_{\alpha\in C}G_{\alpha}$. By Lemma \ref{grouprestrictlemma} $S$ is a Rothberger bounded subset of $\prod_{\alpha<\kappa}G_{\alpha}$. Evidently $S$ generates a Rothberger bounded group.
$\Box$

\begin{corollary}\label{khfailure}
For uncountable cardinals $\kappa$, ${\sf BC}(^{\kappa}2,\kappa)$ implies the failure of ${\sf KH}(\kappa,\aleph_1)$.
\end{corollary}

Solovay proved that Kurepa's Hypothesis is consistent (it holds in the Constructible Universe L). Silver \cite{silver} proved that the negation of Kurepa's Hypothesis is consistent if, and only if, it is consistent that there is an inaccessible cardinal. Thus, the consistency of ${\sf BC}_{\aleph_1}$ (or even ${\sf BC}(^{\omega_1}2, \aleph_1)$) requires the consistency of the existence of an inaccessible cardinal.

\begin{theorem}\label{second} For an uncountable cardinal $\kappa$, each of the statements below implies all the succeeding ones; moreover, if {\sf BC} holds, then they are all equivalent.
\begin{enumerate}
  \item{${\sf BC}_{\kappa}$.}
  \item{${\sf BC}(^{\kappa}2,\kappa)$.} 
  \item{Each Rothberger bounded subgroup of the group $(^{\kappa}2,\oplus)$ has cardinality at most $\kappa$.}
  \item{$\neg{\sf KH}(\kappa,\aleph_1)$.} 
\end{enumerate}
\end{theorem}
{\flushleft{\bf Proof:}} It is clear that (1) implies (2) and that (2) implies (3). Theorem \ref{first} shows that (3) implies (4).  
To see that (4) implies (1), assume that ${\sf BC}_{\aleph_0}$ holds and ${\sf BC}_{\kappa}$ fails. Let $G$ be an $\aleph_0$-bounded group of weight $\kappa$ witnessing the failure of ${\sf BC}_{\kappa}$, and let $S\subseteq G$ be a Rothberger bounded subset of $G$ of cardinality $\kappa^+$. By Theorem \ref{guranquantify} choose separable metrizable groups $G_{\alpha}$, $\alpha<\kappa$ such that $G$ embeds as a topological group into $\prod_{\alpha<\kappa}G_{\alpha}$. Then $S$, considered a subset of $\prod_{\alpha<\kappa}G_{\alpha}$, is Rothberger bounded in the latter and of cardinality $\kappa^+$. Thus, by Lemma \ref{grouprestrictlemma}, for each countable set $C\subseteq \kappa$ the set $S_C\subseteq \prod_{\alpha\in C}G_{\alpha}$ is Rothberger bounded. Since $\prod_{\alpha\in C}G_{\alpha}$ is a separable metrizable space and ${\sf BC}_{\aleph_0}$ holds, Theorem \ref{BCequivs} implies that $S_C$ is countable. Considering $S$ as a family of subsets of $\bigcup S$ we find that $S$ is a witness that the statement ${\sf KH}(\kappa,\aleph_1)$ is true.  
$\Box$

\begin{corollary}\label{BCandKH} Assume ${\sf BC}_{\aleph_0}$. Then ${\sf BC}(^{\omega_1}2,\aleph_1)$ is equivalent to the failure of the Kurepa Hypothesis.
\end{corollary}

A family $\mathcal{F}$ of subsets of an uncountable cardinal $\kappa$ is said to be a $\kappa$-\emph{Kurepa family} if $\vert\mathcal{F}\vert > \kappa$ and for each infinite $\alpha<\kappa$ we have $\vert\{X\cap\alpha:X\in\mathcal{F}\}\vert\le \vert \alpha\vert$. Kurepa's Hypothesis for $\kappa$, ${\sf KH}_{\kappa}$, states that there exists a $\kappa$-Kurepa family. Note that a $(\kappa,\kappa)$-Kurepa family is also an example of a $\kappa$-Kurepa family. Thus, ${\sf KH}(\kappa,\kappa)$ implies ${\sf KH}_{\kappa}$\footnote{We don't know if the converse is true.}.  It is clear that ${\sf KH}(\kappa,\kappa)$ also implies ${\sf KH}(\kappa,\lambda)$ for each uncountable $\lambda<\kappa$.

\begin{lemma}\label{Lemma 12.5} Let $\kappa$ be an uncountable cardinal. If ${\sf KH}(\lambda, \lambda)$ fails for each uncountable $\lambda \le \kappa$, then ${\sf KH}(\kappa, \aleph_1)$ fails.
\end{lemma}

\begin{theorem}\label{third} 
 For an uncountable cardinal $\kappa$, each of the statements below implies all the succeeding ones. If {\sf BC} holds, then statements (1)-(3) are equivalent. If ${\sf BC}_{\lambda}$ holds for each infinite cardinal $\lambda < \kappa$, then all four statements are equivalent:
\begin{enumerate}
  \item{${\sf BC}_{\kappa}$.}
  \item{${\sf BC}(^{\kappa}2,\kappa)$.} 
  \item{$\neg{\sf KH}(\kappa,\aleph_1)$.}   
  \item{$\neg{\sf KH}(\kappa,\kappa)$.}
\end{enumerate}
\end{theorem}
{\flushleft{\bf Proof:}} In light of Theorem \ref{second} and the definitions, the only statement that requires proof is: For each uncountable cardinal $\kappa$, if for each infinite $\lambda<\kappa$, ${\sf BC}_{\lambda}$ holds, then (4) implies (3).  This will be proven by induction on $\kappa$. 

For $\kappa=\aleph_1$ there is nothing to prove. Thus, assume that $\kappa>\aleph_1$ and that the statement has been proven for all uncountable cardinals less than $\kappa$. Towards proving the contrapositive, assume that ${\sf KH}(\kappa,\aleph_1)$ holds. Let $\mathcal{F}$ be a family of subsets of $\kappa$ witnessing ${\sf KH}(\kappa,\aleph_1)$. Then $\vert\mathcal{F}\vert>\kappa$, and for each countable subset $A$ of $\kappa$, $\vert\{X\cap A:X\in\mathcal{F}\}\vert\le\aleph_0$.  Since ${\sf BC}_{\lambda}$ holds for each infinite cardinal $\lambda<\kappa$, Theorem \ref{second} implies that  ${\sf KH}(\lambda,\aleph_1)$ fails for each uncountable $\lambda<\kappa$. By the induction hypothesis, ${\sf KH}(\lambda,\lambda)$ fails for each uncountable cardinal $\lambda<\kappa$. Thus, for each uncountable ordinal $\alpha<\kappa$ we have $\vert\{X\cap\alpha: X\in\mathcal{F}\}\vert\le \vert\alpha\vert$. But this means ${\sf KH}(\kappa,\kappa)$ holds.
$\Box$

Now consider inaccessible cardinals of uncountable cofinality. An uncountable regular cardinal $\kappa$ is \emph{ineffable} if there is for each sequence $(A_{\alpha}:\alpha<\kappa)$ where for each $\alpha$, $A_{\alpha}\subseteq \alpha$, a set $A\subseteq \kappa$ for which $\{\alpha<\kappa: A_{\alpha} = A\cap\alpha\}$ is stationary. 

\begin{theorem}\label{ineffable} Let $\kappa$ be an ineffable cardinal. If ${\sf BC}_{\lambda}$ holds for each infinite cardinal $\lambda<\kappa$, then ${\sf BC}_{\kappa}$ holds.
\end{theorem}
{\bf Proof:} Let $\kappa$ be an ineffable cardinal. Then ${\sf KH}(\kappa,\kappa)$ fails (\cite{KDevlin}, Theorem VII.3.1). By Theorem \ref{third} ${\sf BC}_{\kappa}$ holds.
$\Box$ 

An increasing sequence $(\nu_{\alpha}:\alpha<\mu)$ of cardinals is said to be continuous if for each limit ordinal $\beta<\mu$ we have $\nu_{\beta}=\sup\{\nu_{\alpha}:\alpha<\beta\}$.

\begin{theorem}\label{sslimitunc} Let $\kappa$ be a singular strong limit cardinal of uncountable cofinality $\mu$. If there is an increasing continuous $\mu$-sequence of cardinal numbers $(\nu_{\alpha}:\alpha<\mu)$ below $\kappa$ with supremum equal to $\kappa$ such that $\{\alpha<\mu:\, {\sf BC}_{\nu_{\alpha}}\}$ is a stationary subset of $\mu$, then ${\sf BC}_{\kappa}$.
\end{theorem}
{\bf Proof:} Let $\kappa$ be a singular strong limit cardinal of uncountable cofinality $\mu$. Let $S$ be a Rothberger bounded subset of an $\aleph_0$-bounded group $G$ of weight $\kappa$. By Theorem \ref{guranquantify} we may assume that $G$ embeds as topological group in the product $\Pi_{\alpha<\kappa}G_{\alpha}$ where each $G_{\alpha}$ is a separable metrizable group. Let  $(\nu_{\alpha}:\alpha<\mu)$, an increasing continuous $\mu$-sequence of cardinal numbers with supremum equal to $\kappa$, be such that $\{\alpha<\mu:\, {\sf BC}_{\nu_{\alpha}} \mbox{ holds}\}$ is a stationary subset of $\mu$.

For each $\alpha<\mu$ the set $S_{\alpha} = \{f\lceil_{\nu_{\alpha}}:f\in S\}$ is Rothberger bounded in the $\aleph_0$-bounded group $\prod_{\beta<\nu_{\alpha}}G_{\beta}$ of weight at most $\nu_{\alpha}$. 
By hypothesis the set $\{\alpha<\mu:\vert S_{\alpha}\vert\le\nu_{\alpha}\}$ is stationary. 

Theorem 6 of \cite{EHM} implies that the pairwise disjoint family $(S_{\alpha}:\alpha<\mu)$ has at most $\kappa$ almost disjoint transversals. Since distinct elements of $S$ specify distinct almost disjoint transversals of $(S_{\alpha}:\alpha<\mu\}$, it follows that $\vert S\vert \le\kappa$. $\Box$

Next we explore the relevance of Chang's Conjecture to instances of ${\sf BC}_{\kappa}$. Consider a countable language ${\sf L}$ with a distinguished unary relation symbol ${\sf U}$. We say that a structure $\mathfrak{A}$ of ${\sf L}$ is of type $(\kappa,\lambda)$ if the underlying set $A$ of $\mathfrak{A}$ has cardinality $\kappa$, and $\{x\in A:{\sf U}^{\mathfrak{A}}(x)\}$ has cardinality $\lambda$.

For infinite cardinal numbers $\kappa$, $\lambda$, $\mu$ and $\nu$ the symbol
\begin{equation}\label{genericchang}
(\kappa,\lambda)\changarrow(\mu,\nu)
\end{equation}
denotes the statement that for each countable language ${\sf L}$ with a distinguished unary relation symbol ${\sf U}$, and for each structure $\mathfrak{A}$ of type $(\kappa,\lambda)$ there is an elementary substructure $\mathfrak{B}$ of type $(\mu,\nu)$. 
The instances of interest have $\kappa>\lambda$, $\mu>\nu$, $\kappa\ge\mu$ and $\lambda>\nu$. The instance $(\aleph_2,\aleph_1)\changarrow(\aleph_1,\aleph_0)$ is the classical conjecture of Chang.   

Rowbottom \cite{Rowbottom} discovered a convenient combinatorial equivalent for (\ref{genericchang}): For infinite cardinal numbers $\kappa$, $\lambda$, $\mu$ and $\nu$ the symbol
\begin{equation}\label{rowbottomrelation}
\kappa \rightarrow\lbrack \mu]^{<\aleph_0}_{\lambda,\nu}
\end{equation}
denotes the statement that for each function $f$ from $\lbrack\kappa\rbrack^{<\aleph_0}$, the set of finite subsets of $\kappa$, into $\lambda$, there is a set $X\subseteq \kappa$ such that $\vert X\vert = \mu$, and $\vert\{f(Y): Y \mbox{ is a finite subset of } X\}\vert\le \nu$. The following lemma, a special case of a theorem of Rowbottom, is stated in the form we will use.

\begin{lemma}[Rowbottom]\label{rowbottomth} Let $\kappa>\lambda$ be infinite cardinal numbers. Then  $(\kappa^+,\kappa)\changarrow(\lambda^+,\lambda)$ is equivalent to  $\kappa^+ \rightarrow\lbrack \lambda^+ \rbrack^{<\aleph_0}_{\kappa,\lambda}$.
\end{lemma}

For infinite cardinal numbers $\kappa$, $\lambda$, $\mu$ and $\nu$ the symbol
\begin{equation}\label{rowbottomweakrelation}
\kappa \rightarrow\lbrack \mu]^{2}_{\lambda,\nu}
\end{equation}
denotes the statement that for each function $f$ from $\lbrack\kappa\rbrack^{2}$, the set of 2-element subsets of $\kappa$, into $\lambda$, there is a set $X\subseteq \kappa$ such that $\vert X\vert = \mu$, and $\vert\{f(Y): Y \subseteq X,\, \vert Y\vert=2\}\vert\le \nu$. 

It is evident that $\kappa \rightarrow\lbrack \mu]^{<\aleph_0}_{\lambda,\nu}$ implies $\kappa \rightarrow\lbrack \mu]^{2}_{\lambda,\nu}$. When $\lambda$ is a regular cardinal the converse is also true. A proof of this fact can be gleaned from the corresponding argument for $\kappa=\aleph_2$, $\lambda=\mu=\aleph_1$ and $\nu=\aleph_0$ on page 592 of \cite{BT}.

\begin{theorem}\label{CCandBC} Assume that for the infinite cardinal numbers $\kappa$ and $\lambda$ the partition relation
$\kappa^+ \rightarrow\lbrack \lambda^+ \rbrack^{2}_{\kappa,\lambda}$ holds.
Then ${\sf BC}_{\lambda}$ implies ${\sf BC}_{\kappa}$, and  ${\sf BC}(^{\lambda}2,\lambda)$ implies ${\sf BC}(^{\kappa}2,\kappa)$.
 \end{theorem}
{\flushleft{\bf Proof:}}  
Suppose, towards deriving a contradiction, that  ${\sf BC}_{\kappa}$ fails. 
 Select an $\aleph_0$-bounded group $(G,*)$ of weight $\kappa$ and a subset $X$ of $G$  
such that $X$ is Rothberger bounded and $\vert X\vert = \kappa^+$.  By Theorem \ref{guranquantify} there are separable metrizable groups $(G_{\alpha}:\alpha<\kappa)$ such that $G$ is a subgroup of $\prod_{\alpha<\kappa}G_{\alpha}$, and $X$ is a subset of this product. 
Define a coloring $\Phi$ from $\lbrack X\rbrack^{2}$ 
 to $\kappa$ 
 so that 
\[
  \Phi(\{f,g\}) =\min \{\gamma<\kappa:\,f(\gamma)\neq g(\gamma)\}.
\] 
Apply the partition relation to this coloring to find a subset $Y$ of $X$ and a subset $C$ of $\kappa$ 
 such that $\vert Y\vert = \lambda^+$ and $\vert C\vert = \lambda$ 
 and $\Phi$ restricted to $\lbrack Y\rbrack^{2}$ has values all in $C$. $Y_{C}$ is Rothberger bounded since the projection map is a continuous homomorphism, and $\vert Y_{C}\vert = \lambda^+$ since $\Phi$ is one-to-one on $Y$. 
  But then the  group $G_{C}$ 
   contains a $\lambda^+$-sized Rothberger bounded set $Y_{C}$, and as $G_{C}\subseteq \prod_{\alpha\in C}G_{\alpha}$ this $\aleph_0$-bounded group has weight at most $\lambda$. This provides a contradiction to {\sf BC}$_{\lambda}$. 

The proof that ${\sf BC}(^{\lambda}2,\lambda)$ implies ${\sf BC}(^{\kappa}2,\kappa)$ is left to the reader. $\Box$

\section{Consistency results}

We now consider the consistency of instances of the general Borel Conjecture. 

\subsection{Consistency of the total failure of the general Borel Conjecture}

\begin{lemma}\label{boundedpreserve} If $(G,*)$ is an $\aleph_0$-bounded (totally bounded) topological group and $({\mathbb P},<)$ is a forcing notion, then 
\[
  {\mathbf 1}_{\mathbb P}\forces``(\check{G},*) \mbox{ is $\aleph_0$-bounded (respectively totally bounded)}".
\]
\end{lemma}
{\flushleft{\bf Proof:}} 
Note that the notion of being $\aleph_0$-bounded or of being totally bounded is upwards absolute. 
$\Box$

\begin{theorem}\label{cohenreals}
If $(G,*)$ is an $\aleph_0$-bounded group then in generic extensions by uncountably many Cohen reals, $(G,*)$ is Rothberger bounded.
\end{theorem}
{\flushleft{\bf Proof:}}
Let $(\poset(\kappa),<)$ denote the partially ordered set for adding $\kappa>\aleph_0$ Cohen reals.  By Lemma \ref{boundedpreserve} ${\mathbf 1}_{{\poset}(\kappa)}\forces``(\check{G},*) \mbox{ is $\aleph_0$-bounded}"$. Let $(\dot{\mathcal{U}}_n:n<\omega)$ be a name for a sequence of elements of $\onbd$. Since $\poset(\kappa)$ has the countable chain condition and $\kappa$ is uncountable, there is a countable subset $C$ of $\kappa$ such that $(\dot{\mathcal{U}}_n:n<\omega)$ is a $\poset(C)$-name. As the forcing factors over $C$ we may assume that in fact the sequence so named is a ground model sequence. Since $(G,*)$ is $\aleph_0$-bounded in this model also, we may select for each $n$ a countable set $X_n\subset G$ such that $G=X_n*U_n$. For each $x$ define for each $n$, $f_x(n) = m$ if $x\in x_m*U_n$, $x_m\in X_n$. These objects are all in the ground model. Take a Cohen real over the ground model. It selects a sequence of elements of $G$ which witness Rothberger boundedness.
$\Box$

\begin{theorem}\label{failurecon} It is consistent, relative to the consistency of {\sf ZFC}, that ${\sf BC}(^{\kappa}2,\kappa)$ fails for each infinite cardinal number $\kappa$.
\end{theorem} 
{\flushleft{\bf Proof:}} In the model of Theorem \ref{cohenreals}, for each infinite cardinal $\kappa$ the ground model version of the additive group $^{\kappa}2$ is a Rothberger bounded group of cardinality $2^{\kappa}$. $\Box$

Since adding $\aleph_1$ Cohen reals leaves large cardinal properties of the ground model intact, there is no large cardinal property that implies any instance of ${\sf BC}(^{\kappa}2,\kappa)$.

\subsection{Consistency of ${\sf BC}_{\aleph_0}+{\sf BC}_{\aleph_1}$.}

A partially ordered set $(\poset,<)$ is said to have the \emph{Laver property} if for each $h\in\,^{\omega}\omega$ it is forced that whenever $\tau$, a term in the forcing language of $\poset$, is such that 
$(\forall n)(\tau(n)<\check{h}(n))$ then there exists an $f\in \,^{\omega}(\lbrack\omega\rbrack^{<\omega})$ such that $(\forall n)(\vert f(n)\vert\le 2^n) \mbox{ and }(\forall n)(\tau(n)\in\check{f}(n)))$, and for all but finitely many $n$, $f(n)\subseteq h(n)$.

If in the generic extensions obtained from a partially ordered set all the real numbers are members of the ground model, then the partially ordered set has the Laver property by default. 

The importance of the Laver property is twofold:
\begin{lemma}[Shelah]\label{shelahpreservelaver} A countable support iteration of partially ordered sets, each satisfying the Laver property, satisfies the Laver property\footnote{See Conclusion 2.12 in Chapter VI.2 of \cite{SSPFA}.}.
\end{lemma}
The second important fact about the Laver property is the following folklore result for which a proof can be found in \cite{BJ}, Lemma 3.1:
\begin{lemma}\label{nonsmzpreserve}
Let $X$ be a set of real numbers which does not have strong measure zero. If $(\poset,<)$ is a partially ordered set with the Laver property, then 
\[
  {\mathbf 1}_{\poset}\forces``\check{X}\mbox{ does not have strong measure zero.}"
\]
\end{lemma}

A cardinal $\kappa$ is said to be 1-\emph{inaccessible} if it is inaccessible, and there are $\kappa$ many inaccessible cardinal numbers less than $\kappa$. Now we obtain the following consistency result:
\begin{theorem}\label{conbelowomega2} If it is consistent that there is a  1-inaccessible cardinal, then it is consistent that {\sf ZFC} plus Borel's Conjecture plus the negation of Kurepa's Hypothesis, plus $2^{\aleph_1} = \aleph_2$ hold. 
\end{theorem}
{\flushleft{\bf Proof:}} Let $(\kappa_{\alpha}:\alpha<\kappa)$ be a monotonic enumeration of the inaccessible cardinals below $\kappa$. We construct a $\kappa$-stage countable support iteration $\poset_{\kappa}$ as follows: 
Let $\pi:\kappa\rightarrow\kappa\times\kappa$ be a bookkeeping function such that
\begin{itemize}
  \item{For each $(\beta,\gamma)\in\kappa\times\kappa$ the set $\{\alpha<\kappa:\pi(\alpha)=(\beta,\gamma)\}$ is cofinal in $\kappa$;}
  \item{If $\pi(\alpha)=(\beta,\gamma)$ then $\beta\le\alpha$.} 
\end{itemize}
\underline{$\poset_1$ is defined as follows:}\\
 Let ${\mathbb L}_0$ denote the Levy collapse of $\kappa_0$ to $\omega_2$ with countable conditions. By Silver's Theorem, 
\begin{equation} \mathbf{1}_{{\mathbb L}_{0}}\forces ``{\sf CH } + \mbox{ There are no $\omega_1$ Kurepa trees } "
\end{equation}
Let $\langle\dot{T}^0_{\gamma}:\gamma<\check{\kappa}_0\rangle$ enumerate ${\mathbb L}_0$-names of all $\omega_1$-trees with nodes members of $\omega_1$. By Silver's Theorem each has at most $\aleph_1$ cofinal branches. Pick $\pi(0) = (0,\gamma_0)$, and let $\dot{\mathbb E}_0$ be an ${\mathbb L}_0$-name for a proper partially ordered set that does not add reals and specializes\footnote{In the sense of Baumgartner - see Section 8 of \cite{baumgartner}.} $\dot{T}^0_{\gamma_0}$ (see Chapter 5, Theorem 6.1 and Theorem 7.1 of \cite{SSPFA}). Here we use the fact that if an $\omega_1$-tree has $\le \aleph_1$ cofinal branches of length $\omega_1$, then it has a subtree with no cofinal $\omega_1$-branches, such that rendering this subtree special ensures that no further forcing that preserves $\omega_1$ will add new cofinal $\omega_1$-branches through the original tree. Since ${\mathbf 1}_{{\mathbb L}_0}\forces``\dot{\mathbb E}_0 \mbox{ has the Laver property}"$ it follows that ${\mathbb L}_0*\dot{\mathbb E}_0$ has the Laver property. Next, let $\dot{\mathbb M}$ be a ${\mathbb L}_0*\dot{\mathbb E}_0$-name for the Mathias reals partially ordered set. Since the Mathias reals partially ordered set has the Laver property and forces that every uncountable ground-model set of reals does not have strong measure zero, we find that
${\mathbb L}_0*\dot{\mathbb E}_0*\dot{\mathbb{M}}$ has the Laver property and forces that {\sf CH} holds and every uncountable set of reals from its ground model fails to be strong measure zero.
We set $\poset_1 =  Q_0 = {\mathbb L}_0*\dot{\mathbb E}_0*\dot{\mathbb{M}}$.

With $\alpha\le\kappa$, and assuming that each $\poset_{\beta}$ has been defined for $\beta<\alpha$.\\ 
\underline{$\poset_{\alpha}$ is defined as follows:}\\
\underline{$\alpha=\beta+1$ and $\beta\ge 1$:} Define a $\poset_{\beta}$ name $\dot{Q}_{\beta}$ for a partially ordered set as follows: Let $\dot{\mathbb L}_{\beta}$ be a $\poset_{\beta}$ name for the Levy collapse of $\kappa_{\beta}$ to $\omega_2$ with countable conditions. Let $\langle \dot{T}^{\beta}_{\gamma}:\gamma<\check{\kappa}_{\beta}\rangle$ enumerate $\dot{\mathbb L}_{\beta}$-names for all $\omega_1$ trees with nodes elements of $\omega_1$ With $\pi(\beta) = (\delta,\gamma)$ let $\dot{\mathbb E}_{\beta}$ be a $\dot{\mathbb L}_{\beta}$-name for specializing the $\omega_1$ tree $\dot{T}^{\delta}_{\gamma}$ (note that as $\delta\le\beta$ the most recent Levy collapse ensures that this tree is not a Kurepa tree), and let $\dot{\mathbb M}$ be a $\dot{\mathbb L}_{\beta}*\dot{\mathbb E}_{\beta}$-name for the Mathias reals partially ordered set\footnote{Instead of the Mathias reals partially ordered set, one could also use the Laver reals partially ordered set introduced in \cite{laver}.}. Finally we set
\[
  \dot{Q}_{\beta} = \dot{\mathbb L}_{\beta}*\dot{\mathbb E}_{\beta}*\dot{\mathbb M} \mbox{ and }\poset_{\alpha} = \poset_{\beta}*\dot{Q}_{\beta}.
\]
Then we have 
\begin{equation}
  { \mathbf{1}_{{\mathbb P}_{\beta}}\forces``\dot{Q}_{\beta} \mbox{ has the Laver property}"}
\end{equation}
and also
\begin{equation} 
  \mathbf{1}_{{\mathbb P}_{\beta}}\forces``\dot{\mathbf 1}_{\dot{Q}_{\beta}}\forces \mbox{uncountable ground model sets of reals are not strong measure zero} "
\end{equation}
\underline{$\alpha$ a limit ordinal:} If $\alpha$ has countable cofinality then $\poset_{\alpha}$ is the inverse limit of the $\poset_{\beta}$, $\beta<\alpha$, and else $\poset_{\alpha}$ is the direct limit of $\poset_{\beta},\, \beta<\alpha$.\\

Since $\kappa$ is inaccessible, for each $\beta<\kappa$ $\vert\poset_{\beta}\vert<\kappa$. Then $\poset_{\kappa}$ has the $\kappa$-chain condition. It also follows from Lemma \ref{shelahpreservelaver} that $\poset_{\beta}$, $\beta\le\kappa$ has the Laver property. 

To see that 
\[
  {\mathbf 1}_{\poset_{\kappa}}\forces ``\mbox{ There are no $\omega_1$ Kurepa trees}"
\]
let 
\[
  {\mathbf 1}_{\poset_{\kappa}}\forces``(\check{\omega}_1,\dot{\prec}) \mbox{ is a tree order}"
\]
Since $\poset_{\kappa}$ has the $\kappa$-chain condition and $\kappa$ is inaccessible we find a $\beta<\kappa$ such that 
$(\check{\omega}_1,\dot{\prec})$ is a $\poset_{\beta}$ name and $ {\mathbf 1}_{\poset_{\beta}}\forces``(\check{\omega}_1,\dot{\prec}) \mbox{ is a tree order}"$. But then 
\[
 {\mathbf 1}_{\poset_{\beta}}\forces ``\dot{\mathbf 1}_{\dot{\mathbb L}_{\beta}}\forces ``(\check{\omega}_1,\dot{\prec}) \mbox{ is not a Kurepa tree}""
\]
Now let $\dot{T}^{\beta}_{\gamma}$ be the ${\mathbb L}_{\beta}$ name for $(\check{\omega}_1,\dot{\prec})$, and choose an $\alpha\ge \beta$ such that $\pi(\alpha)=(\beta,\gamma)$, and now consider $\poset_{\alpha+1}$. Since $\dot{\mathbb L}_{\alpha}$ is  a Levy collapse of an inaccessible cardinal Silver's Theorem implies that $(\check{\omega}_1,\dot{\prec})$ is an $\omega_1$ tree with no more than $\aleph_1$ cofinal $\omega_1$ branches. Since $\pi(\alpha)=(\beta,\gamma)$, it follows that $\dot{\mathbb E}_{\alpha}$ specializes $(\check{\omega}_1,\dot{\prec})$. Consequently, 
\[
  {\mathbf 1}_{\poset_{\kappa}}\forces``(\check{\omega}_1,\dot{\prec}) \mbox{ is not a Kurepa tree}"
\]
To see that 
\[
  {\mathbf 1}_{\poset_{\kappa}}\forces ``\mbox{{\sf BC}}"
\]
let $\dot{X}$ be a $\poset_{\kappa}$ name such that
\[
  {\mathbf 1}_{\poset_{\kappa}}\forces ``\dot{X} \mbox{ is a set of real numbers of cardinality }\aleph_1"
\]
By the $\kappa$ chain condition and the strong inaccessibility of $\kappa$ choose a $\beta<\kappa$ such that 
$\dot{X}$ is a $\poset_{\beta}$ name and
\[
  {\mathbf 1}_{\poset_{\beta}}\forces ``\dot{X} \mbox{ is a set of real numbers of cardinality }\aleph_1"
\]
From the definition of $\dot{Q}_{\beta}$ it is clear that
\[
   {\mathbf 1}_{\poset_{\beta}}\forces ``\dot{\mathbf 1}_{\dot{Q}_{\beta}}\forces``\dot{X} \mbox{ is not strong measure zero}""
\]
and thus
\[
   {\mathbf 1}_{\poset_{\beta+1}}\forces ``\dot{X} \mbox{ is not strong measure zero}"
\]
Since $\poset_{\lbrack\beta+2,\kappa)}$ has the Laver property it follows that 
\[
   {\mathbf 1}_{\poset_{\kappa}}\forces ``\dot{X} \mbox{ is not strong measure zero}."
\]
We leave to the reader the standard argument that in the generic extension we have $2^{\aleph_0} = \aleph_2=\kappa$
$\Box$

\begin{corollary}\label{BCtwice} If it is consistent that there is a 1-inaccessible cardinal, then ${\sf BC}_{\aleph_0} + {\sf BC}_{\aleph_1}$ is consistent. 
\end{corollary}
{\flushleft{\bf Proof:}} Corollary \ref{BCandKH} and Theorem \ref{conbelowomega2}. $\Box$

Since we may assume the ground model is ${\mathbf L}$, we may assume that the generic model of Theorem \ref{conbelowomega2} satisfies: For each uncountable cardinal $\kappa$, $2^{\kappa}=\kappa^+$ holds. It is well-known that $2^{\aleph_0} = \aleph_1$ implies $\neg {\sf BC}_{\aleph_0}$. Theorem \ref{conbelowomega2} shows that $2^{\aleph_1}=\aleph_2$ does not imply  $\neg {\sf BC}_{\aleph_1}$.

\subsection{Consistency of $(\forall n<\omega){\sf BC}_{\aleph_n}$.}

\begin{lemma}\label{closednessKHpreservation} Let $\kappa$ and $\lambda$ be uncountable cardinal numbers with $\lambda<\kappa$. 
Let $(\poset,<)$ be a partially ordered set which is $\kappa^+$-closed. If $\neg{\sf KH}_{\lambda}$, then ${\mathbf 1}_{\poset}\forces ``\neg{\sf KH}_{\check{\lambda}}.$" 
\end{lemma}
{\flushleft{\bf Proof:}} This follows from Theorem VII.6.14 of \cite{Kunen}. $\Box$ 

\begin{lemma}\label{closednessBCpreservation} Let $\kappa$ be a regular cardinal number with $2^{\aleph_0}\le\kappa$. Let $(\poset,<)$ be a partially ordered set which is $\kappa$-closed. If ${\sf BC}_{\aleph_0}$, then
${\mathbf 1}_{\poset}\forces ``{\sf BC}_{\aleph_0}.$" 
\end{lemma}
{\flushleft{\bf Proof:}} This also follows from Theorem VII.6.14 of \cite{Kunen}: No new sets of real numbers of cardinality $\aleph_1$ are added by this forcing.  $\Box$

\begin{theorem}\label{fragment} Assume it is consistent that the following three statements hold: ${\sf BC}_{\aleph_0}$,  $\neg{\sf KH}_{\aleph_1}$, $2^{\aleph_1}=\aleph_2$, and there are inaccessible cardinals $\kappa_0<\cdots<\kappa_n<\cdots$, $n<\omega$. Then ${\sf BC}_{\aleph_0} + (\forall n<\omega)(0<n\Rightarrow \neg{\sf KH}_{\aleph_n})$ is consistent. 
\end{theorem}
{\flushleft{\bf Proof:}} As in Exercise (F4) on p. 295 of \cite{Kunen} define a countable support iterated forcing poset $\poset$ such that successively ``for each $n[>0]$, $\kappa_n$ is Levy collapsed to $\aleph_{n+2}$ by conditions of cardinality $\le \kappa_{n-1}$". Then, by \cite{baumgartner}, Theorem 2.5,  
$(\poset,<)$ is $\aleph_2$-closed. By Lemmas \ref{closednessKHpreservation} and \ref{closednessBCpreservation} $(\poset,<)$  preserves ${\sf BC}_{\aleph_0}\, +\, \neg{\sf KH}_{\aleph_1}$. By the cited exercise from \cite{Kunen}, in the resulting generic extension we have $(\forall n<\omega)(0<n\Rightarrow \neg{\sf KH}_{\aleph_n})$. $\Box$

\begin{corollary}\label{fragment2}
If ${\sf BC}_{\aleph_0}$ $+$ $\neg{\sf KH}_{\aleph_1}$ $+$ $2^{\aleph_1}=\aleph_2$ $+$ there are inaccessible cardinals $\kappa_0<\cdots<\kappa_n<\cdots$, $n<\omega$ is consistent, then $(\forall n<\omega){\sf BC}_{\aleph_n}$ is consistent. 
\end{corollary}
{\flushleft{\bf Proof:}}  Theorem \ref{fragment} and Theorem \ref{third}. $\Box$

\subsection{Consistency of ${\sf BC}_{\kappa}$ first failing at $\kappa=\aleph_{\omega}$.}

\begin{theorem}\label{breaksatalephomega} Suppose it is consistent that there is an inaccessible cardinal $\kappa$ such that there are $\kappa$ inaccessible cardinals below $\kappa$, and $\omega$ inaccessible cardinals above $\kappa$. Then it is consistent that $(\forall n<\omega){\sf BC}_{\aleph_n}$ while also $\neg {\sf BC}_{\aleph_{\omega}}$.
\end{theorem}
{\flushleft{\bf Proof:}} We may assume the ground model is ${\mathbf{L}}$. Let $\lambda$ be the limit of the inaccessible cardinals assumed to exist in the hypothesis. Thus, $\lambda$ has countable cofinality and there is a Kurepa family on $\lambda$. Performing the forcing in Theorem \ref{conbelowomega2}, followed by the forcing in Theorem \ref{fragment} preserves this Kurepa family, but collapses $\lambda$ to $\aleph_{\omega}$. $\Box$

\subsection{Consistency of ${\sf BC}_{\aleph_{\omega}}$.}
\begin{quote}
\end{quote}

An uncountable cardinal number $\kappa$ is said to be $\mu$-strong if there is an elementary embedding $j:\, {\sf V}\rightarrow {\sf M}$ with critical point $\kappa$ such that ${\sf V}_{\mu}\subseteq {\sf M}$. $\kappa$ is said to be a \emph{strong} cardinal if it is $\mu$-strong for all $\mu$. 
\begin{theorem}\label{BCandstrongcards} If it is consistent that for an uncountable cardinal $\kappa$ of countable cofinality both $2^{\kappa}=\kappa^+$ and ${\sf BC}(^{\kappa}2,\kappa)$, then it is consistent that there is a strong cardinal. 
\end{theorem}
{\flushleft{\bf Proof:}}
Todorcevic proved (see for example Chapter 7 of \cite{stevo}) that if $\kappa$ is an uncountable cardinal of countable cofinality then $\Box_{\kappa}$ 
 plus ${\sf cof}(\lbrack\kappa\rbrack^{\aleph_0},\subseteq) = \kappa^+$ implies that there is a cofinal in $\lbrack\kappa\rbrack^{\aleph_0}$ family of countable sets that witnesses ${\sf KH}(\kappa,\aleph_1)$. Applying Theorem \ref{second} we find that $\neg {\sf BC}(^{\kappa}2,\kappa)$ holds. 

Thus, if ${\sf BC}(^{\kappa}2,\kappa)$ and $2^{\kappa}=\kappa^+$ hold, then $\Box_{\kappa}$ fails. Jensen has proved that failure of $\Box_{\kappa}$ for uncountable $\kappa$ of countable cofinality implies the existence of an inner model with a strong cardinal (see Fact 2.6 of \cite{CFM}). 
$\Box$

In consistency strength strong cardinals lie between measurable cardinals and strongly compact cardinals: A strong cardinal is measurable. If $\kappa$ is strongly compact then $\Box_{\lambda}$ fails for each cardinal $\lambda>\kappa$ and thus there is an inner model with a strong cardinal.

\begin{theorem}\label{BCandstrongcards2} If ${\sf BC}(^{\kappa}2, \kappa)$ holds for an uncountable cardinal $\kappa$ of countable cofinality for which we have $\lambda^{\aleph_0}<\kappa$ for all $\lambda<\kappa$, then the axiom of projective determinacy is true. 
\end{theorem}
{\flushleft{\bf Proof:}} Let $\kappa$ be an uncountable cardinal of countable cofinality such that for each cardinal $\lambda<\kappa$ we have $\lambda^{\aleph_0}<\kappa$. Also assume that ${\sf BC}(^{\kappa}2,\kappa)$ holds. Using the argument in the proof of Theorem \ref{BCandstrongcards}, it follows that $\Box_{\kappa}$ fails. But this, by \cite{SZ} Corollary 6, implies that the axiom of projective determinacy is true. $\Box$

\begin{corollary}\label{alephomegaBCPD} If $2^{\aleph_0}<\aleph_{\omega}$ and if ${\sf BC}(^{\aleph_{\omega}}2,\aleph_{\omega})$, 
then Projective Determinacy holds.
\end{corollary}

\begin{corollary}\label{alephomegaBCLAD} If for each $n<\omega$ we have $2^{\aleph_n}<\aleph_{\omega}$ and if ${\sf BC}(^{\aleph_{\omega}}2,\aleph_{\omega})$, then Determinacy holds in {\sf L}$(\reals)$.
\end{corollary}
{\flushleft{\bf Proof:}} The argument is as in the proof of Theorem \ref{BCandstrongcards2}, except that we now use \cite{Steel}, Theorem 0.1, which states that if there is a singular strong limit cardinal $\kappa$ such that $\Box_{\kappa}$ fails, then the axiom of determinacy holds in ${\sf L}(\reals)$. $\Box$

Now we determine upper bounds on the consistency strength of ${\sf BC}_{\aleph_{\omega}}$. 

\begin{lemma}\label{closednesspreserves} Let $\kappa>\lambda$ be infinite cardinal numbers. Let $(\poset,<)$ be a $\kappa^{++}$-closed partially ordered set. If the partition relation $\kappa^+ \rightarrow\lbrack \lambda^+ \rbrack^{2}_{\kappa,\lambda}$ holds, then
\begin{equation}\label{closedpreservefla}
  {\mathbf 1}_{\poset}\forces``\check{\kappa}^+ \rightarrow\lbrack \check{\lambda}^+ \rbrack^{2}_{\check{\kappa},\check{\lambda}}.
\end{equation}
\end{lemma}
{\flushleft{\bf Proof:}} Assume that (\ref{closedpreservefla}) is false. Choose $p\in\poset$ such that
$p\forces ``\check{\kappa}^+ \not\rightarrow\lbrack \check{\lambda}^+ \rbrack^{2}_{\check{\kappa},\check{\lambda}}"$.
Choose a $\poset$-name $\dot{f}$ such that
 $p\forces ``\dot{f}:\lbrack\check{\kappa}^+\rbrack^{2}\rightarrow\check{\lambda}^+ \mbox{ witnesses } \check{\kappa}^+ \not\rightarrow\lbrack \check{\lambda}^+ \rbrack^{2}_{\check{\kappa},\check{\lambda}}"$.
Since $\poset$ is $\kappa^{++}$ closed we find a ground model function $g:\lbrack\kappa^+\rbrack^{2}\rightarrow \lambda^+$ and a $q<p$ such that 
  $q\forces ``\dot{f}=\check{g} \mbox{ on the ground model set }\lbrack\kappa^+\rbrack^{2}"$.
Applying the partition relation $\kappa^+ \rightarrow\lbrack \lambda^+ \rbrack^{2}_{\kappa,\lambda}$ in the ground model to $g$ we find an uncountable ground model subset $S$ of $\kappa^+$ and a ground model subset $C\subset \lambda^+$ such that $\vert C\vert\le\lambda$ and for any $F\in \lbrack S\rbrack^2$ we have $g(F)\in C$. But then
  $q\forces ``\mbox{For each finite }F\in \lbrack\check{S}\rbrack,\, \dot{f}(F)\in\check{C}"$,
 contradicting the fact that $q<p$ and $p\forces ``\check{\kappa}^+ \not\rightarrow\lbrack \check{\lambda}^+ \rbrack^{2}_{\check{\kappa},\check{\lambda}}"$.
$\Box$

For an ordinal number $\alpha$ and for a cardinal number $\kappa$, the notation $\kappa^{+\alpha}$ denotes the $\alpha$-th cardinal number larger than $\kappa$.

\begin{lemma}[Levinski, Magidor, Shelah]\label{GeneralCollapse}
Let $\mu$ be an inaccessible cardinal. Let $(\poset,<)$ be a partially ordered set such that $\vert \poset\vert = \mu$ and $\poset$ has the $\mu$-chain condition. If $(\mu^{+\omega+1},\mu^{+\omega})\changarrow(\aleph_1,\aleph_0)$ holds, then
${\mathbf 1}_{\poset}\forces``(\check{\mu}^{+\omega+1},\check{\mu}^{+\omega})\changarrow(\check{\aleph}_1,\check{\aleph}_0)"$.
\end{lemma}
{\flushleft{\bf Proof:}} See page 168 of \cite{LMS}. $\Box$ 

If $\mu$ is an inaccessible cardinal and $\theta$ is a regular cardinal with $\mu>\theta$, then
${\sf Lv}(\mu,\theta)$ is the set of $p$ such that $p$ is a function with ${\sf dom}(p)\subseteq \mu\times \theta$,  $\vert p\vert<\theta$, and for all $(\alpha,\xi)\in {\sf dom}(p)$ we have $p(\alpha,\xi)\in \alpha$. For $p$ and $q$ in ${\sf Lv}(\mu,\theta)$ write $q<p$ if $p\subset q$. Then the partially ordered set $({\sf Lv}(\mu,\theta),<)$ is the \emph{Levy collapse}. It is $\theta$-closed, has the $\mu$-chain condition, and has cardinality $\mu$.

\begin{lemma}[Levinski, Magidor, Shelah]\label{finLevyCollapse} Let $\lambda>\mu$ be inaccessible cardinals. If 
\[
  (\lambda^{+\omega+1},\lambda^{+\omega})\changarrow(\mu^{+\omega+1},\mu^{+\omega})
\]
 holds, then ${\mathbf 1}_{({\sf Lv}(\mu^{+\omega},\omega),<)}\forces``(\check{\lambda}^{+\omega+1},\check{\lambda}^{+\omega})\changarrow({\aleph}_1,{\aleph}_0)"$.
\end{lemma}
{\flushleft{\bf Proof:}} See p. 168 of \cite{LMS}. $\Box$

\begin{lemma}\label{LevyCollapse}
Let $\mu$ be an inaccessible cardinal such that $(\mu^{+\omega+1},\mu^{+\omega})\changarrow(\aleph_1,\aleph_0)$ holds. Let $\alpha<\mu$ be an ordinal and let $\theta<\mu$ be the cardinal $\aleph_{\omega\cdot\alpha+2}$. Then
\[
  {\mathbf 1}_{({\sf Lv}(\mu,\theta),<)}\forces``({\aleph}_{\omega\cdot(\alpha+1)+1},{\aleph}_{\omega\cdot(\alpha+1)})\changarrow(\check{\aleph}_1,\check{\aleph}_0)"
\]
\end{lemma}
{\flushleft{\bf Proof:}} By Lemma \ref{GeneralCollapse}, ${\mathbf 1}_{({\sf Lv}(\mu,\theta),<)}\forces``(\check{\mu}^{+\omega+1},\check{\mu}^{+\omega})\changarrow(\check{\aleph}_1,\check{\aleph}_{0})"$. It is well-known that ${\mathbf 1}_{({\sf Lv}(\mu,\theta),<)}\forces``\vert\check{\mu}\vert={\aleph}_{\omega\cdot\alpha+3}"$. Consequently, ${\mathbf 1}_{({\sf Lv}(\mu,\theta),<)}\forces``\vert\check{\mu}^{+\omega+1}\vert={\aleph}_{\omega\cdot\alpha+\omega+1}"$ and ${\mathbf 1}_{({\sf Lv}(\mu,\theta),<)}\forces``\vert\check{\mu}^{+\omega}\vert={\aleph}_{\omega\cdot\alpha+\omega}"$
$\Box$
\vspace{0.1in}

For $0<n<\omega$ the uncountable cardinal $\lambda$ is said to be $n$-\emph{huge} if there 
 is an elementary embedding $j:V\rightarrow M$ to a transitive inner model $M$ of the set theoretic universe $V$ such that $\lambda$ is the critical point of $j$, and setting $\kappa_0=\lambda$ and $\kappa_{i+1}=j(\kappa_i)$ for $i< n$, we have $^{\kappa_n}M\subseteq M$. It can be shown that if $\lambda$ is $n$-huge then, in the above notation, each of the cardinals $\kappa_i$, $0\le i\le n$ is measurable. 
\begin{lemma}[Levinski, Magidor, Shelah]\label{hugeCC} If $\lambda$ is a 2-huge cardinal, then 
\[
  (\kappa_1^{+\omega+1},\kappa_1^{+\omega})\changarrow(\lambda^{+\omega+1},\lambda^{+\omega}).
\]
\end{lemma}

\begin{theorem}\label{alephomandBC} If it is consistent that there is a 2-huge cardinal, then it is consistent that ${\sf BC}_{\aleph_{\omega}}$.
\end{theorem}
{\flushleft{\bf Proof: }} 
Let $\lambda$ be a 2-huge cardinal and let $j$ be an elementary embedding witnessing this. Put $\kappa = j(\lambda)$. It is known that $\kappa$ is measurable and $\kappa>\lambda$. By Lemma \ref{hugeCC}
$(\kappa^{+\omega+1},\kappa^{+\omega})\changarrow(\lambda^{+\omega+1},\lambda^{+\omega})$.

Forcing first with $\poset_0 = ({\sf Lv}(\lambda^{+\omega},\omega),<)$ we obtain by Lemma \ref{finLevyCollapse} a generic extension in which we have 
$({\kappa}^{+\omega+1},{\kappa}^{+\omega})\changarrow({\aleph}_1,{\aleph}_0)$. 
Since $\kappa$ is still measurable in this generic extension, it is an inaccessible limit of inaccessible cardinals. Now let $\poset_1$ be the corresponding partially ordered set for Theorem \ref{conbelowomega2}. Then $\vert\poset_1\vert=\kappa$, and $\poset_1$ has the $\kappa$-chain condition. 

Since ${\mathbf 1}_{\poset_1}\forces``\check{\kappa}=\aleph_2"$, Lemma \ref{GeneralCollapse} gives $(\aleph_{\omega+1},\aleph_{\omega})\changarrow(\aleph_1,\aleph_0)$ in the generic extension. By Theorem \ref{conbelowomega2} this generic extension also satisfies ${\sf BC}_{\aleph_0}$ and ${\sf BC}_{\aleph_1}$. But then since ${\sf BC}_{\aleph_0}$ holds, Lemma \ref{rowbottomth} and Theorem \ref{CCandBC} imply that ${\sf BC}_{\aleph_{\omega}}$ holds in this generic extension.
$\Box$

By the facts in the table below ${\sf BC}_{\aleph_{\omega\cdot n}} \, +\, {\sf BC}_{\aleph_{\omega\cdot n +1}}$ holds  in the model of Theorem \ref{alephomandBC} for each $n<\omega$.
\begin{center}
\begin{tabular}{|l|l|} \hline
       The argument to prove Lemma \ref{hugeCC} gives for all n:    & Then the model of Theorem \ref{alephomandBC} gives: \\ \hline
      & \\
      $(\kappa_1^{+\omega+n+1},\kappa_1^{+\omega+n})\changarrow(\lambda^{+\omega+n+1},\lambda^{+\omega+n}).$ & $(\aleph_{\omega+n+1},\aleph_{\omega+n})\changarrow(\aleph_{n+1},\aleph_{n}).$\\
       & \\
      $(\kappa_1^{+\omega (n+1)+1},\kappa_1^{+\omega (n+1)})\changarrow(\lambda^{+\omega (n+1)+1},\lambda^{+\omega (n+1)}).$ & $(\aleph_{\omega(n+1)+1},\aleph_{\omega(n+1)})\changarrow(\aleph_{\omega n+1},\aleph_{\omega n}).$\\ \hline
\end{tabular}
\end{center}

\subsection{Consistency of: For a proper class of $\kappa$ with $cf(\kappa)=\aleph_0$, ${\sf BC}_{\kappa}$.}
\begin{quote}\end{quote}
For ordinal number $\alpha$ a cardinal number $\lambda$ is said to be $n$-huge $\alpha$ times if there is for each ordinal $\beta<\alpha$ an elementary embedding $j_{\beta}$ into a transitive inner model $M_{\beta}$ such that each $j_{\beta}$ witnesses that $\lambda$ is $n$-huge, and when $\beta<\delta<\alpha$, then $j_{\beta}(\lambda)<j_{\delta}(\lambda)$.

\begin{theorem}[Barbanel, Di Prisco, Tan]\label{n+1tonhuge} If $\lambda$ is $n+1$-huge, then there is a cardinal $\mu<\lambda$ such that $\mu$ is $n$-huge, and 
\[
  \{\alpha<\lambda: \mbox{ There is an $n$-huge elementary embedding $j$ with }j(\mu)=\alpha\}
\]
is a stationary subset of $\lambda$.
\end{theorem}

\begin{corollary}\label{manyCC}
If $\lambda$ is a 3-huge cardinal, then there is a 2-huge cardinal $\mu$ such that
\[
 T =  \{\alpha<\lambda: \alpha \mbox{ is measurable and }(\alpha^{+\omega+1},\alpha^{+\omega})\changarrow(\mu^{+\omega+1},\mu^{+\omega})\}
\]
is a stationary subset of $\lambda$.
\end{corollary}
{\flushleft{\bf Proof:}} Lemma \ref{hugeCC} and Theorem \ref{n+1tonhuge}. $\Box$

\begin{corollary}\label{properclassBC}
If it is consistent that there is a 3-huge cardinal, then it is consistent that ${\sf BC}_{\aleph_0}$ as well as ${\sf BC}_{\aleph_1}$, and there is a proper class of cardinals $\kappa$ such that $\omega={\sf cf}(\kappa)$, and ${\sf BC}_{\kappa}$ as well as ${\sf BC}_{\kappa^+}$. 
\end{corollary}
{\flushleft{\bf Proof:}}
Now let $T$ be as in Corollary \ref{manyCC}. Upon forcing with $({\sf Lv}(\mu^{+\omega},\omega),<)$ we find that 
\[
 T =  \{\alpha<\lambda: \alpha \mbox{ is a measurable cardinal and }(\alpha^{+\omega+1},\alpha^{+\omega})\changarrow(\aleph_1,\aleph_0)\}
\]
Enumerate $T$ in increasing order as $(\alpha_{\xi}:\xi<\lambda)$. Next we force with the poset of Theorem \ref{conbelowomega2}, using an iteration of length $\alpha_0$. In the resulting model we have $\alpha_0 = \aleph_2$ and for all $\xi>0$, $\alpha_{\xi}$ is still measurable. Moreover we have for each $\xi$ that ${\sf BC}_{\alpha_{\xi}^{+\omega}}$ as well as ${\sf BC}_{\alpha_{\xi}^{+\omega+1}}$ hold. Since $\lambda$ is still measurable, $V_{\lambda}$ is a model of ZFC, and in $V_{\lambda}$ we have for each $0<\xi<\lambda$ that ${\sf BC}_{\alpha_{\xi}^{+\omega}}$ as well as ${\sf BC}_{\alpha_{\xi}^{+\omega+1}}$ hold.
$\Box$

\section{Questions} 

In Theorem \ref{failurecon} we showed that $(\forall \kappa)(\neg {\sf  BC}(^{\kappa}2,\kappa))$ holds in generic extensions by $\aleph_1$ Cohen reals. 
\begin{problem} Does ${\mathbf V}={\mathbf L}$ imply $(\forall \kappa)(\neg {\sf BC}(^{\kappa}2,\kappa)$?
\end{problem}

In all our models of instances of ${\sf BC}_{\kappa}$ also ${\sf BC}_{\aleph_0}$ is true.

\begin{problem} Is it consistent that ${\sf BC}(^{\kappa}2, \kappa)$ holds for some uncountable cardinal $\kappa$, while ${\sf BC}$ fails? What if $\kappa = \aleph_1$ or $\kappa = \aleph_{\omega}$?
\end{problem}

${\sf BC}(^{\kappa}2,\kappa)$ implies that every Rothberger subspace of $^{\kappa}2$ has cardinality at most $\kappa$. For $\kappa=\aleph_0$ the converse is true. This is not known for $\kappa>\aleph_0$.

\begin{problem} Is it for each infinite cardinal $\kappa$ true that if each Rothberger subspace of $^{\kappa}2$ has cardinality at most $\kappa$, then ${\sf BC}(^{\kappa}2,\kappa)$ holds?
\end{problem}

For each $\kappa$, ${\sf BC}_{\kappa}$ implies ${\sf BC}(^{\kappa}2,\kappa)$. For $\kappa>\aleph_0$ it is not clear if the converse is true.

\begin{problem} Is it true that for each uncountable cardinal $\kappa$, ${\sf BC}(^{\kappa}2,\kappa)$ implies ${\sf BC}_{\kappa}$?
\end{problem}

We obtained from the consistency of a large cardinal hypothesis the consistency of the statement that ${\sf BC}_{\kappa}$ holds for a proper class of cardinals $\kappa$ (of countable cofinality). 

\begin{problem}
Is ${\sf ZFC}+(\forall \kappa){\sf BC}(^{\kappa}2,\kappa)$  consistent relative to the consistency of any large cardinal axioms?
\end{problem}

Our findings indicate that ${\sf BC}(^{\aleph_\omega}2,\aleph_{\omega})$ has considerable consistency strength.

\begin{problem} What is the exact consistency strength of ${\sf BC}(^{\aleph_{\omega}}2,\aleph_{\omega})$?
\end{problem}

\section{Acknowledgments}

We thank Justin Moore for bringing to our attention the fact that there are forcing notions which specialize $\omega_1$ trees  and do not add reals. We also thank Frank Tall for discussions on various methods for collapsing cardinals.


\begin{thebibliography}{}

\bibitem{BDT} J.B. Barbanel, C.A. Di Prisco and I.B. Tan, \emph{Many-times huge and superhuge cardinals}, {\bf The Journal of Symbolic Logic} 49:1 (1984), 112 - 122 

\bibitem{BJ} T. Bartoszynski and H. Judah, \emph{Strong measure zero sets}, {\bf Israel Mathematical Conference Proceedings} 6 (1993), 13 - 62 

\bibitem{baumgartner} J.E. Baumgartner, \emph{Iterated Forcing}, in Surveys in Set Theory {\bf London Mathematical Society Lecture Notes} 87 (1983), 1 - 59.

\bibitem{BT} J.E. Baumgartner and A.D. Taylor, \emph{Saturation properties of ideals in generic extensions {\tt II}}, {\bf Transactions of the American Mathematical Society} 271:2 (1982), 587 - 609.

\bibitem{Borel} E. Borel, \emph{Sur la classification des ensembles de mesure nulle}, {\bf Bulletin de la Societe Mathematique de France} 47 (1919), 97 - 125. 

\bibitem{CarlsonBC} T.J. Carlson, \emph{Strong measure zero and strongly meager sets}, {\bf Proceedings of the American Mathematical Society} 118:2 (1993), 577 - 586.

\bibitem{CFM} J. Cummings, M. Foreman and M.Magidor, \emph{Squares, scales and stationary reflection}, {\bf Journal of Mathematical Logic} 1:1 (2001), 35 - 98. 

\bibitem{KDevlin} K.J. Devlin, \emph{Constructibility}, {\bf Springer-Verlag}, 1984. 

\bibitem{EHM} P. Erd\H{o}s, A. Hajnal and E.C. Milner, \emph{On sets of almost disjoint subsets of a set}, {\bf Acta Mathematica Academiae Scientiarum Hungaricae} 19 (1968), 209-218

\bibitem{Guran} I.I. Guran, \emph{On topological groups close to being Lindel\"of}, {\bf Soviet Math. Dokl.} 23 (1981), 173 - 175.

\bibitem{HS} A. Halko and S. Shelah, \emph{On strong measure zero subsets of $^{\kappa}2$}, {\bf Fundamenta Mathematicae} 170:3 (2001), 219 - 229.

\bibitem{coc2} W. Just, A.W. Miller, M. Scheepers and P.J. Szeptycki, \emph{The Combinatorics of Open Covers (II)}, {\bf Topology and its Applications} 73 (1996), 241 - 266. 

\bibitem{Kelley} J.L. Kelley, \emph{General Topology}, {\bf Springer-Verlag} Graduate Texts in Mathematics 27, 1975

\bibitem{Kunen} K.Kunen, \emph{Set Theory: An introduction to independence proofs}, {\bf Studies in Logic and Mathematics} 102 1980.

\bibitem{laver} R. Laver, \emph{On the consistency of Borel's conjecture}, {\bf Acta Mathematica} 137 (1976), 151 - 169.

\bibitem{LMS} J.-P. Levinski, M. Magidor and S. Shelah, \emph{Chang's conjecture for $\aleph_{\omega}$}, {\bf Israel Journal of Mathematics} 69:2 (1990), 161 - 172.

\bibitem{rothberger38} F. Rothberger, \emph{Eine Versch\"arfung der Eigenschaft C}, {\bf Fundamenta Mathematicae} 30 (1938), 50 - 55. 

\bibitem{Rowbottom} F. Rowbottom, \emph{Large cardinals and small constructible sets}, {\bf Annals of Mathematical Logic} 3 (1971), 1-44.

\bibitem{SZ} E. Schimmerling and M. Zeman, \emph{Square in core models}, {\bf Bulletin of Symbolic Logic} 7:3 (2001), 305 - 314.

\bibitem{SSPFA} S. Shelah, \emph{Proper Forcing}, {\bf Springer-Verlag: Lecture Notes in Mathematics} 940, 1982.

\bibitem{sierpinski} W. Sierpi\'nski, \emph{Sur un ensemble non d\'enombrable, dont toute image continue est de mesure nulle}, {\bf Fundamenta Mathematicae} 11 (1928), 301--304.

\bibitem{silver} J.H. Silver, \emph{The independence of Kurepa's conjecture and two-cardinal conjectures in model theory}, {\bf AMS Proc. Symp. Pure Math.} 13:1 (1971), 383 - 390.

\bibitem{Steel} J.R. Steel, \emph{{\sf PFA} implies ${\sf AD}^{{\sf L}(\reals)}$}, {\bf The Journal of Symbolic Logic} 70:4 (2005), 1255 - 1296.

\bibitem{SM} E. Szpilrajn-Marczewski, \emph{La dimension et la mesure}, {\bf Fundamenta Mathematicae} 28 (1937), 81-89

\bibitem{MT} M. Tkachenko, \emph{Introduction to topological groups}, {\bf Topology and its Applications} 86 (1998), 179 - 231. 

\bibitem{stevo} S. Todorcevic, \emph{Walks on ordinals and their characteristics}, {\bf Birkh\"auser Verlag series: Progress in Mathematics} 263, 2007

\bibitem{TW} B. Tsaban and T. Weiss, \emph{Products of special sets of real numbers}, {\bf Real Analysis Exchange} 30 (2004/5), 819 - 836.

\end{thebibliography}
\end{document}